\documentstyle[12pt]{amsart}
\def\opn#1#2{\def#1{\operatorname{#2}}} 
\opn\chara{char}
\opn\length{\ell}
\opn\pd{pd}
\opn\rk{rk}
\opn\projdim{proj\,dim}
\opn\injdim{inj\,dim}
\opn\rank{rank}
\opn\depth{depth}
\opn\grade{grade}
\opn\height{height}
\opn\embdim{emb\,dim}
\opn\codim{codim}

\opn\Tr{Tr}
\opn\bigrank{big\,rank}
\opn\superheight{superheight}\opn\lcm{lcm}
\opn\trdeg{tr\,deg}%
\opn\reg{reg}
\opn\lreg{lreg}
\opn\skel{skel}
\opn\com{com}
\opn\div{div}
\opn\Div{Div}
\opn\cl{cl}
\opn\Cl{Cl}
\opn\Spec{Spec}
\opn\Supp{Supp}
\opn\supp{supp}
\opn\Sing{Sing}
\opn\Ass{Ass}
\opn\Ann{Ann}
\opn\Rad{Rad}
\opn\Soc{Soc}
\opn\Ker{Ker}
\opn\Coker{Coker}
\opn\Im{Im}
\opn\Hom{Hom}
\opn\Tor{Tor}
\opn\Ext{Ext}
\opn\End{End}
\opn\Aut{Aut}
\opn\id{id}

\opn\nat{nat}
\opn\pff{pf}
\opn\Pf{Pf}
\opn\GL{GL}
\opn\SL{SL}
\opn\mod{mod}
\opn\ord{ord}
\opn\aff{aff}
\opn\con{conv}
\opn\relint{relint}
\opn\st{st}
\opn\lk{lk}
\opn\cn{cn}
\opn\core{core}
\opn\vol{vol}
\opn\link{link}
\opn\star{star}
\opn\gr{gr}

\def\pot#1#2{#1[\kern-0.28ex[#2]\kern-0.28ex]}

\opn\dirlim{\underrightarrow{\lim}}
\opn\inivlim{\underleftarrow{\lim}}
\def\Implies{\ifmmode\Longrightarrow \else
     \unskip${}\Longrightarrow{}$\ignorespaces\fi}
\def\implies{\ifmmode\Rightarrow \else
     \unskip${}\Rightarrow{}$\ignorespaces\fi}
\def\iff{\ifmmode\Longleftrightarrow \else
     \unskip${}\Longleftrightarrow{}$\ignorespaces\fi}

\let\:=\colon
\newtheorem{Theorem}{Theorem}[section]
\newtheorem{Lemma}[Theorem]{Lemma}

\newtheorem{Remark}[Theorem]{Remark}

\let\epsilon\varepsilon
\let\phi=\varphi
\let\kappa=\varkappa
\def\qed{\ifhmode\textqed\fi
   \ifmmode\ifinner\quad\qedsymbol\else\dispqed\fi\fi}
\def\textqed{\unskip\nobreak\penalty50
    \hskip2em\hbox{}\nobreak\hfil\qedsymbol
    \parfillskip=0pt \finalhyphendemerits=0}
\def\dispqed{\rlap{\qquad\qedsymbol}}

\def\FF{{\cal F}}

\opn\initial{in}
\opn\inim{inm}
\opn\rev{rev}
\opn\Gin{Gin}
\opn\Lex{Lex}
\opn\Shift{Shift}
\opn\shift{shift}
\opn\rate{rate}
\opn\Mon{Mon}
\opn\lex{lex}
\opn\rev{rev}
\opn\red{red}
\opn\max{max}
\opn\min{min}
\opn\initial{in}
\opn\Ker{Ker}
\opn\GL{GL}
\opn\proj{proj}
\begin{document}
\title{
Kruskal--Katona type theorems for clique complexes
arising from chordal and strongly chordal graphs
}
\author{J\"urgen Herzog, Takayuki Hibi, Satoshi Murai,\\
Ng\^o Vi\^et Trung and Xinxian Zheng}
\address{J\"urgen Herzog, Fachbereich Mathematik und
Informatik,
Universit\"at Duisburg--Essen, Campus Essen,
45117 Essen, Germany}
\email{juergen.herzog@@uni-essen.de}
\address{Takayuki Hibi,
Department of Pure and Applied Mathematics,
Graduate School of Information Science and Technology,
Osaka University,
Toyonaka, Osaka 560-0043, Japan}
\email{hibi@@math.sci.osaka-u.ac.jp}
\address{Satoshi Murai,
Department of Pure and Applied Mathematics,
Graduate School of Information Science and Technology,
Osaka University,
Toyonaka, Osaka 560-0043, Japan}
\email{s-murai@@ist.osaka-u.ac.jp}
\address{Ng\^o Vi\^et Trung,
Institute of Mathematics, Vien Toan Hoc,
18 Hoang Quoc Viet, 10307 Hanoi}
\email{nvtrung@@math.ac.vn}
\address{Xinxian Zheng, Fachbereich Mathematik und
Informatik,
Universit\"at Duisburg--Essen, Campus Essen,
45117 Essen, Germany}
\email{xinxian.zheng@@uni-essen.de}
\maketitle
\begin{abstract}
A forest is the clique complex of
a strongly chordal graph and a quasi-forest
is the clique complex of a chordal graph.
Kruskal--Katona type theorems for forests, quasi-forests,
pure forests and pure quasi-forests
will be presented.
\end{abstract}
\section*{Introduction}
Recently, in commutative algebra,
the forest (\cite{F})
and the quasi-forest (\cite{Z} and \cite{HHZ})
have been extensively studied.  Each of these concepts
is, however, well known in combinatorics
(\cite{MM}).
In fact, a forest is the clique complex of
a strongly chordal graph
and a quasi-forest
is the clique complex of a chordal graph.
(A chordal graph is a finite graph
for which every cycle of length $ > 3$ has a chord.
A strongly chordal graph is a chordal graph
for which every cycle of even length $\geq 6$
has a chord that joins two vertices of the cycle
with an odd distance $ > 1$ in the cycle.)

Besides the celebrated $g$-conjecture for spheres
(\cite[pp. 75--76]{Stanley}),
one of the most important open problems
in the study of $f$-vectors of simplicial complexes
is the classification of $f$-vectors of flag
complexes.  (A flag complex is the clique complex
of a finite graph.)
Works in the reserach topic include
\cite{CharneyDavis},
\cite{FerrFroe},
\cite{Frohmader},
\cite{Gal} and \cite{Renteln}.
On the other hand,
the study of $f$-vectors of clique complexes of
chordal graphs was done in
\cite{Eckhoff},
\cite{Kalai1} and
\cite{Kalai2}.

The purpose of the present paper is to give
a Kruskal--Katona type theorem for forests and
quasi-forests (Theorem \ref{Osaka})
as well as
a Kruskal--Katona type theorem
for pure forests and pure quasi-forests
(Theorem \ref{Stockholm}).
These theorems will be proved in Section $2$.
We then show in Section $3$ that
the $f$-vector of a pure quasi-forest
is unimodal.

\section{Kruskal--Katona type theorems}
Let $[n] = \{ 1, \ldots, n \}$ be the vertex set
and $\Delta$ a simplicial complex on $[n]$.
Thus $\Delta$ is a collection of subsets of $[n]$
with the properties that
(i) $\{ i \} \in \Delta$ for each $i \in [n]$ and
(ii) if $F \in \Delta$ and $G \subset F$,
then $G \in \Delta$.
Each element $F \in \Delta$ is a face of $\Delta$.
Let $d = \max \{ |F| : F \in \Delta \}$, where
$|F|$ is the cardinality of $F$.  Then
$\dim \Delta$, the dimension of $\Delta$,
is $d - 1$.  A facet is a maximal face of $\Delta$ under inclusion.
We write $\FF(\Delta)$ for the set of facets of
$\Delta$.
A simplicial complex is called pure if all facets
have the same cardinality.
Let $f_i$ denote the number of faces
$F$ with $|F| = i + 1$.  The vector
$f(\Delta) = (f_0, f_1, \ldots, f_{d-1})$
is called the $f$-vector of $\Delta$.
In particular $f_0 = n$.
If $\{ F_{i_1}, \ldots, F_{i_q} \}$ is
a subset of $\FF(\Delta)$, then we write
$\langle F_{i_1}, \ldots, F_{i_q} \rangle$
for the subcomplex of $\Delta$
whose faces are those faces $F$ of $\Delta$
with $F \subset F_{i_j}$ for some $1 \leq j \leq q$.

A facet $F$ of a simplicial complex $\Delta$ is called
a {\em leaf} if there is a facet $G \neq F$ of $\Delta$,
called a {\em branch} of $F$, such that $H \bigcap F
\subset G \bigcap F$ for all facets $H$ of $\Delta$
with $H \neq F$.  A {\em quasi-forest} is a simplicial
complex $\Delta$ which enjoys an ordering
$F_1, F_2, \ldots, F_s$ of the facets of $\Delta$,
called a {\em leaf order},
such that for each $1 < j \leq s$ the facet
$F_j$ is a leaf of the subcomplex
$\langle F_{1}, \ldots, F_{j-1}, F_j \rangle$
of $\Delta$.  A {\em quasi-tree} is a quasi-forest
which is connected.  A {\em forest} is a simplicial
complex $\Delta$ which enjoys the property that
for every subset $\{ F_{i_1}, \ldots, F_{i_q} \}$
of $\FF(\Delta)$ the subcomplex
$\langle F_{i_1}, \ldots, F_{i_q} \rangle$
of $\Delta$ has a leaf.
A {\em tree} is a forest which is connected.

We now come to Kruskal--Katona type theorems for forests,
quasi-forests, pure forests and pure quasi-forests.

\begin{Theorem}
\label{Osaka}
Given a finite sequence $(f_0, f_1, \ldots, f_{d-1})$
of integers with each $f_i > 0$,
the following conditions are equivalent:

(i) there is a quasi-forest $\Delta$ of dimension $d - 1$
with $f(\Delta) = (f_0, f_1, \ldots, f_{d-1})$;

(ii) there is a forest $\Delta$ of dimension $d - 1$
with $f(\Delta) = (f_0, f_1, \ldots, f_{d-1})$;

(iii) the sequence $(c_1, \ldots, c_{d})$
defined by the formula
\begin{eqnarray}
\label{NewYork}
\sum_{i=0}^{d} f_{i-1} (x - 1)^{i}
= \sum_{i=0}^{d} c_i x^i,
\end{eqnarray}
where $f_{-1} = 1$,
satisfies
$\sum_{i=k}^{d} c_i > 0$ for each $1 \leq k \leq d$.
%


(iv) the sequence
$(b_1,\dots,b_d)$ defined by the formula
\begin{eqnarray}
\label{Toyonaka}
\sum_{i=1}^d f_{i-1} (x-1)^{i-1} = \sum_{i=1}^d b_i x^{i-1}
\end{eqnarray}
is positive, i.e., $b_i > 0$ for $1 \leq i \leq d$.
\end{Theorem}

\begin{Theorem}
\label{Stockholm}
Given a finite sequence $(f_0, f_1, \ldots, f_{d-1})$
of integers with each $f_i > 0$,
the following conditions are equivalent:

(i) there is a pure quasi-forest $\Delta$ of dimension $d - 1$
with $f(\Delta) = (f_0, f_1, \ldots, f_{d-1})$;

(ii) there is a pure forest $\Delta$ of dimension $d - 1$
with $f(\Delta) = (f_0, f_1, \ldots, f_{d-1})$;


(iii)
the sequence $(c_1, \ldots, c_{d})$
defined by (\ref{NewYork}) satisfies
$\sum_{i=k}^{d} c_i > 0$ for each $1 \leq k \leq d$
and $c_i \leq 0$ for each $1 \leq i < d$.

(iv) the sequence
$(b_1,\dots,b_d)$ defined by the formula
(\ref{Toyonaka})
satisfies
$0 < b_1 \leq b_2 \leq \cdots \leq b_d$.
\end{Theorem}

In Section $2$, after preparing
Lemmata \ref{quasiforest}, \ref{Copenhagen}
and \ref{forest}, we will prove
both of Theorems \ref{Osaka} and \ref{Stockholm}
simultaneously.

\section{$f$-vectors of forests and quasi-forests}
We begin with

\begin{Lemma}
\label{quasiforest}
Let $\Delta$ be a quasi-forest on $[n]$
with $s + 1$ facets and
$F_{s+1}, F_{s}, \ldots, F_{1}$
its leaf order.
For each $1 \leq j \leq s$
we write $G_{j}$ for a branch
of the leaf $F_{j}$ in the subcomplex
$\langle F_{s+1}, F_{s}, \ldots, F_{j} \rangle$
of $\Delta$.
Let $\delta_j = |F_j|$ and
$e_j = |F_{j} \bigcap G_j|$.
Let $\dim \Delta = d - 1$ and
let $f(\Delta) = (f_0, f_1, \ldots, f_{d-1})$
be the $f$-vector of $\Delta$.

(a)
One has
\begin{eqnarray}
\label{Sydney}
\sum_{i=0}^{d} f_{i-1} x^{i}
= \sum_{j=1}^{s+1} (1 + x)^{\delta_j} -
\sum_{j=1}^{s} (1 + x)^{e_j},
\end{eqnarray}
where $f_{-1} = 1$.

(b)
Let $k_1 \cdots k_s k_{s+1}$ be a permutation
of $[s+1]$ with
$0 < \delta_{k_1} \leq \cdots \leq \delta_{k_{s}}
\leq \delta_{k_{s+1}} = d$ and
$\ell_1 \cdots \ell_s$ a permutation
of $[s]$ with
$0 \leq e_{\ell_1} \leq \cdots \leq e_{\ell_s}$.
Then $e_{\ell_j} < \delta_{k_j}$
for all $1 \leq j \leq s$.
\end{Lemma}

\begin{pf}
(a)
Let $s = 0$.
Then
$\delta_1 = d$, $f_{i - 1} = {d \choose i}$
and
$\sum_{i=0}^{d} f_{i-1} x^{i}
= (1 + x)^{\delta_1}$.
Let $s \geq 1$ and $\Delta' =
\langle F_{s+1}, F_s, \ldots, F_2 \rangle$.
Let $\dim \Delta' = d' - 1$
and $f(\Delta') = (f'_0, f'_1, \ldots, f'_{d'-1})$.
Since $\Delta'$ is a quasi-forest,
working by induction on the number of facets,
it follows that
$\sum_{i=0}^{d'} f'_{i-1} x^{i}
= \sum_{j=2}^{s+1} (1 + x)^{\delta_j} -
\sum_{j=2}^{s} (1 + x)^{e_j}$.
Now, the number of faces $F$ of
$\Delta$ with $F \not\in \Delta'$
and with $|F| = i$ is
${\delta_{1} \choose i} - {e_1 \choose i}$.
Hence
\[
\sum_{i=0}^{d} f_{i-1} x^{i}
= \sum_{i=0}^{d'} f'_{i-1} x^{i}
+ ((1 + x)^{\delta_{1}} - (1 + x)^{e_1}),
\]
as desired.

(b)
Let $k_p = 1$ and $\ell_q = 1$.
Then $e_{\ell_q} < \delta_{k_p}$.
Since
$\langle F_{s+1}, F_s, \ldots, F_2 \rangle$
is a quasi-forest,
working by induction on the number of facets,
it follows that
(i) in case of $p \leq q$, one has
$e_{\ell_j} < \delta_{k_j}$
for each $1 \leq j < p$
and for each $q < j \leq s$,
and
$e_{\ell_j} \leq e_{\ell_q}
< \delta_{k_p} \leq \delta_{k_j}$
for each $p \leq j \leq q$,
and that (ii) in case of $q < p$,
one has
$e_{\ell_j} < \delta_{k_j}$
for each $1 \leq j < q$
and for each $p < j \leq s$, and
$e_{\ell_j} \leq e_{\ell_{j+1}}
< \delta_{k_{j}} \leq \delta_{k_{j+1}}$
for each $q \leq j < p$.
Hence
$e_{\ell_j} < \delta_{k_j}$
for all $1 \leq j \leq s$.
\end{pf}

\begin{Lemma}
\label{Copenhagen}
Given a quasi-forest of dimension $d - 1$
with $s + 1$ facets and with $f$-vector
$f(\Delta) = (f_0, f_1, \ldots, f_{d-1})$,
there exist finite sequences
$(\delta_1, \ldots, \delta_t, \delta_{t+1})$
and $(e_1, \ldots, e_t)$
of integers, where $0 < t \leq s$
and where
$\delta_i \neq e_j$ for all $i$ and $j$,
satisfying
\[
0 < \delta_1 \leq \cdots \leq \delta_t \leq
\delta_{t+1} = d,
\, \, \, \, \, \, \, \, \, \,
0 \leq e_1 \leq \cdots \leq e_t < d
\]
and
\[
e_j < \delta_j, \, \, \, \, \,
1 \leq j \leq t
\]
which enjoys the formula
\[
\sum_{i=0}^{d} f_{i-1} x^{i}
= \sum_{j=1}^{t+1} (1 + x)^{\delta_j} -
\sum_{j=1}^{t} (1 + x)^{e_j},
\]
where $f_{-1} = 1$.
\end{Lemma}

\begin{pf}
In
Lemma \ref{quasiforest} (b),
in case that $\delta_{k_a} = e_{\ell_b}$
for some $a$ and $b$,
one has $a < b$ and
$e_{\ell_j} < \delta_{k_{j+1}}$
for all $a \leq j < b$.
Hence we can replace
$\delta_{k_1}, \ldots, \delta_{k_{s+1}}$
and
$e_{\ell_1}, \ldots, e_{\ell_s}$
with
$\delta_{k_1}, \ldots, \delta_{k_{a-1}},
\delta_{k_{a+1}}, \ldots, \delta_{k_s},
\delta_{k_{s+1}}$
and
$e_{\ell_1}, \ldots, e_{\ell_{b-1}}, e_{\ell_{b+1}},
\ldots, e_{\ell_s}$.
\end{pf}

\begin{Lemma}
\label{forest}
Let
$(\delta_1, \ldots, \delta_s, \delta_{s+1})$
and $(e_1, \ldots, e_s)$
be sequences of integers,
where $s \geq 0$
and where
$\delta_i \neq e_j$ for all $i$ and $j$,
satisfying
\[
0 < \delta_1 \leq \cdots \leq \delta_s \leq
\delta_{s+1} = d,
\, \, \, \, \, \, \, \, \, \,
0 \leq e_1 \leq \cdots \leq e_s < d
\]
and
\[
e_j < \delta_j, \, \, \, \, \,
1 \leq j \leq s.
\]
Then there is a forest $\Delta$ on $[n]$,
where
$n = \sum_{j=1}^{s+1} \delta_j - \sum_{j=1}^{s} e_j$,
of dimension $d - 1$ such that
the $f$-vector
$f(\Delta) = (f_0, f_1, \ldots, f_{d-1})$
of $\Delta$ satisfies (\ref{Sydney}).
\end{Lemma}

\begin{pf}
First, we construct the subsets
$F_1, \ldots, F_s, F_{s+1}$ of $[n]$,
where $|F_j| = \delta_j$ for each $1 \leq j \leq s + 1$
and
where $|F_j \bigcap F_{j+1}| = e_j$
for each $1 \leq j \leq s$.
Let $F_{s+1} = \{n - d + 1, n - d + 2, \ldots, n \}$.
If we obtain the facet
$F_j = \{ q_1, q_2, \ldots, q_{\delta_j - 1},
q_{\delta_j} \}$,
where $1 \leq q_1 < q_2 < \cdots < q_{\delta_j} \leq n$,
then the facet $F_{j-1}$ is defined to be
\[
F_{j-1} = \{ q_1 - (\delta_{j-1} - e_{j-1}),
\ldots, q_1 - 2, q_1 - 1,
q_{\delta_j - e_{j-1} + 1}, \ldots, q_{\delta_j - 1},
q_{\delta_j} \}.
\]
A crucial property of $F_1, \ldots, F_s, F_{s+1}$
is that
\begin{eqnarray}
\label{crucial}
F_j \bigcap F_{s+1} = \cdots =
F_j \bigcap F_{j+2} = F_j \bigcap F_{j+1}
\end{eqnarray}
for each $1 \leq j \leq s$.

Now, write $\Delta$ for the simplicial complex on $[n]$
of dimension $d - 1$
with $\FF(\Delta) = \{ F_1, \ldots, F_s, F_{s+1} \}$.
It follows from the property (\ref{crucial}) that
$\Delta$ is a quasi-forest with
$F_{s+1}, F_s, \ldots, F_1$ its leaf order.
In addition, for each $1 \leq j \leq s$,
in the quasi-forest
$\langle F_{s+1}, F_s, \ldots, F_j \rangle$
the facet $F_{j+1}$ is a branch of the leaf $F_j$.
Since $|F_j| = \delta_j$ for each $1 \leq j \leq s + 1$
and
$|F_j \bigcap F_{j+1}| = e_j$
for each $1 \leq j \leq s$,
Lemma \ref{quasiforest}
guarantees that
the $f$-vector
$f(\Delta) = (f_0, f_1, \ldots, f_{d-1})$
of $\Delta$ satisfies (\ref{Sydney}).

Finally, we claim that $\Delta$ is a forest.
Let $1 \leq j_1 < j_2 < \cdots < j_q \leq s + 1$
and $\Gamma = \langle F_{j_1}, \ldots, F_{j_q} \rangle$.
Then
\[
F_{j_1} \bigcap F_{j_q} = \cdots =
F_{j_1} \bigcap F_{j_3} = F_{j_1} \bigcap F_{j_2}.
\]
Hence $F_{j_1}$
is a leaf of $\Gamma$ with $F_{j_2}$
its branch.
\end{pf}

We now prove
both of Theorems \ref{Osaka} and \ref{Stockholm}
simultaneously.

\begin{pf}
((ii) $\Rightarrow$ (i))
Since a (resp. pure) forest is a
(resp. pure) quasi-forest,
the $f$-vectors of (resp. pure) forests coincide with
the $f$-vectors of (resp. pure) quasi-forests.

((i) $\Rightarrow$ (iii))
Let $\Delta$ be a quasi-forest
of dimension $d - 1$ with $s + 1$ facets
and $f(\Delta) = (f_0, f_1, \ldots, f_{d-1})$
its $f$-vector.
With the same notation as in Lemma \ref{Copenhagen},
it follows that
\[
\sum_{i=0}^{d} c_i x^{i}
= \sum_{j=1}^{t+1} x^{\delta_j} -
\sum_{j=1}^{t} x^{e_j}.
\]
Thus
\[
\sum_{i=k}^{d} c_i
=
|\{ j : \delta_j \geq k \}|
- |\{ j : e_j \geq k \}| > 0
\]
for each $1 \leq k \leq d$, as desired.

If, in addition, $\Delta$ is pure, then
Lemma \ref{quasiforest} guarantees the existence
of a sequence $(e_1, \ldots, e_s)$ of integers
with $0 \leq e_1 \leq \cdots \leq e_s < d$ such that
\[
\sum_{i=0}^{d} c_i x^{i}
= (s + 1) x^{d} -
\sum_{j=1}^{s} x^{e_j}.
\]
Thus
$\sum_{i=k}^{d} c_i > 0$
for each $1 \leq k \leq d$
and $c_i \leq 0$ for each $1 \leq k < d$.

((iii) $\Rightarrow$ (ii))
Let $s + 1 = \sum_{c_i > 0} c_i$.
Since $\sum_{i=0}^{d} c_i = f_{-1} = 1$,
one has $\sum_{c_i < 0} ( - c_i ) = s$.
Since $\sum_{i=1}^{d} c_i > 0$
and $\sum_{i=0}^{d} c_i = f_{-1} = 1$,
one has $c_0 \leq 0$.
Write
\[
\sum_{i=0}^{d} c_i x^{i}
= \sum_{j=1}^{s+1} x^{\delta_j} -
\sum_{j=1}^{s} x^{e_j},
\]
where $e_i \neq \delta_j$ for all $i$ and $j$,
and where
\[
0 < \delta_1 \leq \cdots \leq \delta_s \leq
\delta_{s+1} = d,
\, \, \, \, \, \, \, \, \, \,
0 \leq e_1 \leq \cdots \leq e_s < d.
\]
Suppose $\sum_{i=k}^{d} c_i > 0$
for each $1 \leq k \leq d$.
We claim
$e_j < \delta_j$ for each $1 \leq j \leq s$.
To see why this is true,
let $j$ denote the biggest integer $\leq s$
for which $\delta_j < e_j$.  Then
$\sum_{i=e_j}^{d} c_i = 0$, a contradiction.

Now, it turns out that the sequences
$(\delta_1, \ldots, \delta_s, \delta_{s+1})$
and $(e_1, \ldots, e_s)$ enjoy the properties
required in Lemma \ref{forest}.
Thus there exists a forest $\Delta$
of dimension $d - 1$ whose $f$-vector
$f(\Delta)$ satisfies (\ref{Sydney}).  In other words, the $f$-vector
$f(\Delta)$ must coincide with the given
sequence $(f_0, f_1, \ldots, f_{d-1})$.

If, in addition, $c_i \leq 0$ for each $1 \leq i < d$,
then $\delta_j = d$ for all $1 \leq j \leq s + 1$.
Hence the forest which is constructed in the proof of
Lemma \ref{forest} is pure.

((iii) $\Leftrightarrow$ (iv))
Let $b_k = \sum_{i=k}^d c_i$ for $k=0,1,\dots,d$.
In other words, the sequence
$(b_0,b_1,\dots,b_d)$ is defined by the formula
$$\sum_{i=0}^d c_i x^i = \sum_{i=1}^d b_i x^{i-1} (x-1) + b_0.$$
Since $b_0 = c_0 + \cdots + c_d = f_{-1}$, 
it follows from (\ref{NewYork}) that
$(b_1,\dots,b_d)$ satisfies (\ref{Toyonaka}).
It is now clear that (iii) is equivalent to (iv), as desired.
\end{pf}

\begin{Remark}
\rm
Let $\Delta$ be a simplicial complex
on $[n]$ of dimension $d - 1$ and
$f(\Delta) = (f_0, f_1, \ldots, f_{d-1})$
its $f$-vector.  Recall from
\cite{BruHer}, \cite{Hibi} and \cite{Stanley}
that the $h$-vector
$h(\Delta) = (h_0, h_1, \ldots, h_d)$
of $\Delta$ is defined by the formula
\begin{eqnarray*}
\sum_{i=0}^{d} f_{i-1}(x - 1)^{d-i}
= \sum_{i=0}^{d} h_i x^{d-i},
\end{eqnarray*}
or equivalently, by the formula
\begin{eqnarray}
\label{general}
\sum_{i=0}^dh_ix^i(1+x)^{d-i}=\sum_{i=0}^df_{i-1}x^{i}.
\end{eqnarray}
In particular $h_0 = 1$ and $h_1 = n - d$.
It follows from
(\ref{NewYork}) and
(\ref{general})
that
\begin{eqnarray*}
\sum_{i=0}^d c_i x ^{i}
=\sum_{i=0}^d h_i x^{d-i} (x-1)^i.
\end{eqnarray*}
Hence
$$c_i = \sum_{j=0}^d (-1)^{d-i} { j \choose d-i} h_j.$$
Consequently, for each $k=1,2,\dots,d$, one has   
$$
b_k = 
\sum_{i=k}^d c_i = \sum_{j=0}^d \left \{ \sum_{i=k}^d (-1)^{d-i} { j 
\choose d-i} h_j \right\}
=1+ \sum_{j=1}^d (-1)^{d-k} { j-1 \choose d-k} h_j.
$$
\end{Remark}

\section{Unimodality of $f$-vectors}
A finite sequence $(a_1, a_2, \ldots, a_N)$ of integers
with each $a_i > 0$ is called {\em unimodal}
if $a_0 \leq \cdots \leq a_j \geq a_{j+1} \geq \cdots
\geq a_N$ for some $0 \leq j \leq N$.
We recall the following well-known

\begin{Lemma}
\label{binomial}
Let $d$ and $e$ be integers with $0 \leq e < d$
and define the sequence $(a_0, a_1, \ldots, a_\delta)$
by the formula
$\sum_{i=0}^{\delta} a_i x^i = (1 + x)^d - (1 + x)^e$.
Then
\[
a_0 \leq a_1 \leq \cdots \leq
a_{[(d+1)/2]} \geq a_{[(d+1)/2]+1}
\geq \cdots \geq a_{d}.
\]
\end{Lemma}

\begin{Theorem}
\label{Kyoto}
The $f$-vector $f(\Delta) = (f_0, f_1, \ldots, f_{d-1})$
of a pure quasi-forest $\Delta$ of dimension $d - 1$
is unimodal.
\end{Theorem}

\begin{pf}
Let $\Delta$ be a pure quasi-forest of domension $d - 1$
with $s+1$ facets and
$f(\Delta) = (f_0, f_1, \ldots, f_{d-1})$
its $f$-vector.
Lemma \ref{quasiforest} says that
there is a sequence $(e_1, \ldots, e_s)$ of integers
with $0 \leq e_1 \leq \cdots \leq e_s < d$ such that
\[
\sum_{i=0}^{d} f_{i-1} x^{i}
= (1 + x)^{d} +
\sum_{j=1}^{s} ((1 + x)^{d} - (1 + x)^{e_j}).
\]
By using Lemma \ref{binomial} it follows that
\[
f_{-1} \leq f_0 \leq \cdots \leq
f_{[(d+1)/2]-1} \geq f_{[(d+1)/2]}
\geq \cdots \geq f_{d-1}.
\]
Hence $(f_0, f_1, \ldots, f_{d-1})$
is unimodal.
\end{pf}

The $f$-vector of a pure simplicial complex
is not necessarily unimodal.
In fact, there is a simplicial convex polytope
such that
the $f$-vector of its boundary complex is not
unimodal (\cite{Bjorner}).
On the other hand, it is proved
in \cite{Whatcanbesaid} that
the $f$-vector
$f(\Delta) = (f_0, f_1, \ldots, f_{d-1})$
of a pure simplicial complex $\Delta$
of dimension $d - 1$ satisfies
\[
f_i \leq f_{d - 2 - i}, \, \, \, \, \,
-1 \leq i \leq [d/2]-1,
\]
together with
\[
f_{-1} \leq f_0 \leq \cdots \leq
f_{[d/2]-1}.
\]

\end{document}